\documentclass[11pt]{article}
\usepackage{fullpage}
\usepackage{graphicx}
\usepackage{amssymb}
\usepackage{epstopdf}

\usepackage{fullpage}

\def\qed{\hfill
\ifhmode\unskip\nobreak\fi\quad\ifmmode\Box\else$\Box$\fi\\ }
\newtheorem{theorem}{Theorem}

\newtheorem{defn}[theorem]{Definition}

\newtheorem{claim}[theorem]{Claim}
\newtheorem{fact}[theorem]{Fact}
\title{Ore's Conjecture for $k=4$ and Gr\" otzsch Theorem}
\author{Alexandr  Kostochka\thanks{
University of Illinois at Urbana--Champaign, Urbana, IL 61801, USA and
 Sobolev Institute of Mathematics, Novosibirsk 630090, Russia. Email:
 kostochk@math.uiuc.edu.
Research of this author
is supported in part by NSF grant DMS-0965587 and by
grants 12-01-00448 and 12-01-00631 of the Russian Foundation for Basic Research. }
\and
Matthew Yancey\thanks{Department of Mathematics, University of Illinois, Urbana,
IL 61801, USA. E-mail: yancey1@illinois.edu.
Research of this author is partially supported by the
Arnold O. Beckman Research Award of the University of Illinois
at Urbana-Champaign and from National Science Foundation grant DMS 08-38434 ``EMSW21-MCTP: Research
Experience for Graduate Students.''}}

\begin{document}

\maketitle

\begin{abstract}
A graph $G$ is  $k$-{\em critical} if it has chromatic number $k$, but every
proper subgraph of $G$ is $(k-1)$--colorable.
Let $f_k(n)$ denote the minimum number of edges in an $n$-vertex $k$-critical graph.
In a very recent paper,
we gave a lower bound, $f_k(n)  \geq F(k,n)$, that is sharp for every $n=1\,({\rm mod }\, k-1)$.
It is also sharp for $k=4$ and every $n\geq 6$. In this note, we present a simple
proof of the bound for $k=4$. It implies the case  $k=4$ of the conjecture by Ore from 1967 that
for every $k\geq 4$ and $n\geq k+2$, $f_k(n+k-1)=f(n)+\frac{k-1}{2}(k - \frac{2}{k-1})$.
We also show that our result implies a
 simple short proof of the Gr\" otzsch Theorem that every triangle-free planar graph is $3$-colorable.
\end{abstract}

\section{Introduction}

A {\em proper $k$-coloring}, or simply $k$-{\em coloring}, of a graph $G = (V, E)$ is a function
$f:V \rightarrow \{1,2,\dots,k\}$ such that for each $uv \in E$, $f(u) \neq f(v)$.
A graph $G$ is $k$-{\em colorable} if there exists a $k$-coloring of $G$. The {\em chromatic number}, $\chi(G)$,
of a graph $G$ is the smallest $k$ such that $G$ is $k$-{colorable}.
A graph $G$ is $k$-{\em critical} if $G$ is not $(k-1)$-colorable, but  every proper subgraph of $G$ is $(k-1)$-colorable.
  Then every $k$-critical
graph has chromatic number $k$ and every $k$-chromatic graph contains
a $k$-critical subgraph.


The only $1$-critical graph is $K_1$, and the only $2$-critical graph is $K_2$.
The only $3$-critical graphs are the odd cycles.
Let $f_k(n)$ be the minimum number of edges in a $k$-critical graph with $n$ vertices.
Since $\delta(G)\geq k-1$ for every $k$-critical $n$-vertex graph $G$,
$\quad f_k(n)\geq \frac{k-1}{2}n\quad$
 for all $n\geq k$, $n\neq k+1$. Equality is achieved
for $n=k$ and for $k=3$ and $n$ odd.
In 1957, Dirac~\cite{D1} asked to determine
$f_k(n)$ and proved  that for $k\geq 4$ and
$n\geq k+2$,
$\quad f_k(n)\geq \frac{k-1}{2}n+\frac{k-3}{2}$.
The bound is tight for $n=2k-1$.
Gallai~{\cite{G2}} found exact values of $f_k(n)$ for $k+2\leq n\leq 2k-1$:

\begin{theorem}[Gallai~\cite{G2}] \label{gallai1}
If $k \geq 4$ and $k+2\leq n\leq 2k-1$, then
$$f_k(n)=\frac{1}{2}  \left((k-1)n+(n-k)(2k-n)\right)-1.$$
 \end{theorem}
He also proved that $\quad f_k(n)\geq \frac{k-1}{2}n+\frac{k-3}{2(k^2-3)}n\quad$
for all $k\geq 4$ and
$n\geq k+2$.
Gallai in 1963 and Ore~\cite{O} in 1967 reiterated the question on finding $f_k(n)$.
Ore observed that Haj\'{o}s' construction implies
\begin{equation}\label{upper f_k}
f_k(n + k - 1) \leq f_k(n) + \frac{(k-2)(k+1)}2 = f_k(n) + (k-1)(k - \frac{2}{k-1})/2,
\end{equation}
which yields that $\phi_k:=\lim_{n\to \infty}\frac{f_k(n)}{n}$ exists and satisfies
$\; \phi_k\leq \frac{k}{2}-\frac{1}{k-1}$.
Ore~\cite{O} also conjectured that for every $n\geq k+2$, in
(\ref{upper f_k}) equality holds.

More detail on known results about $f_k(n)$ and Ore's Conjecture the reader can find in~\cite{J}[Problem 5.3]
and our recent paper~\cite{KY}.
In~\cite{KY} we proved the following bound.

\begin{theorem} \label{k-critical}
If $k \geq 4$ and $G$ is $k$-critical, then
$ |E(G)|\geq  \left\lceil \frac{(k+1)(k-2)|V(G)|-k(k-3)}{2(k-1)}\right\rceil$. In other words, if
$k\geq 4$ and $n\geq k,\,n\neq k+1$, then
\begin{equation}\label{j20}
f_k(n)\geq  F(k,n):=\left\lceil \frac{(k+1)(k-2)n-k(k-3)}{2(k-1)}\right\rceil.
\end{equation}
\end{theorem}

This bound is exact for $k=4$ and every $n\geq 6$. For every $k\geq 5$, the bound is exact for every
$n \equiv 1\, ({\rm mod}\ k-1)$, $n\neq 1$. In particular, $\phi_k=\frac{k}{2}-\frac{1}{k-1}$ for every $k\geq 4$.
The result also confirms the above conjecture by Ore from 1967 for $k=4$ and every $n\geq 6$
and also for $k\geq 5$ and all
$n \equiv 1 \, ({\rm mod} \ k-1)$, $n\neq 1$. One of the corollaries of Theorem~\ref{k-critical} is a
   short proof of the following theorem due to Gr\" otzsch~\cite{Gr}:

\begin{theorem}[\cite{Gr}]\label{gr} Every triangle-free planar graph is $3$-colorable.
\end{theorem}

The original proof of Theorem~\ref{gr} is somewhat sophisticated. There were subsequent simpler proofs (see, e.g.~\cite{Tho} and references
therein),
but Theorem~\ref{k-critical} yields a half-page proof.
A disadvantage of this  proof is that
 the proof of Theorem~\ref{k-critical} itself is not too simple. The goal of this note is to give
a simpler proof of the case $k=4$ of Theorem~\ref{k-critical} and to deduce Gr\" otzsch' Theorem from
this result. Note that even the case $k=4$ was a well-known open problem (see, e.g.~\cite{J2}[Problem 12]
and recent paper~\cite{FM}). 
Some further consequences for coloring planar graphs are discussed in~\cite{BKY}.

In Section 2 we prove Case $k=4$ of Theorem~\ref{k-critical} and in Section 3 deduce  Gr\" otzsch Theorem from it.
Our notation is standard. In particular, $\chi(G)$ denotes the chromatic number of graph $G$,
$G[W]$ is the subgraph of a graph  $G$ induced by the vertex set $W$. For
a vertex $v$ in a graph $G$, $d_G(v)$ denotes the degree
of vertex $v$ in graph $G$, $N_G(v)$ is the set of neighbors of $v$. 
If the graph $G$ is clear from the context, we drop the subscript.


\section{Proof of Case $k=4$ of Theorem \ref{k-critical}}
The case $k=4$ of Theorem \ref{k-critical} can be restated as follows.

\begin{theorem} \label{4-critical}
If  $G$ is $4$-critical, then
$ |E(G)|\geq  \left\lceil \frac{5|V(G)|-2}{3}\right\rceil$. \end{theorem}

\begin{defn}
For $R \subseteq V(G)$, define \emph{the potential of $R$} to be
$\quad \rho_{G}(R) = 5|R| - 3|E(G[R])|$.
When there is no chance for confusion, we will use $\rho(R)$.
Let $P(G) = \min_{\emptyset \neq R \subseteq V(G)} \rho(R)$.
\end{defn}

\begin{fact} \label{list of potentials} 
We have
	 $\rho_{K_1}(V(K_1)) = 5$,
	 $\rho_{K_2}(V(K_2)) = 7$,	
	 $\rho_{K_{3}}(V(K_{3})) = 6$,
	 $\rho_{K_4}(V(K_4)) = 2$.
\end{fact}


Note that $|E(G)|\geq\frac{5|E(G)|-2}{3}$ is equivalent to $\rho(V(G)) \leq 2$.
Suppose Theorem~\ref{4-critical} does not hold.
Let $G$ be a vertex-minimal $4$-critical graph with  $\rho(V(G)) > 2$.
This implies that
\begin{equation}\label{j26}
\mbox{if $|V(H)|<|V(G)|$ and $P(H) > 2$, then $H$ is $3$-colorable.}
\end{equation}

\begin{defn}\label{def1} For a graph $G$, a set $R\subset V(G)$ and a $3$-coloring $\phi$ of
$G[R]$, the graph $Y(G,R,\phi)$ is constructed as follows. First, for $i=1,2,3$, let
$R'_i$ denote the set of vertices in $V(G)-R$ adjacent to at least one vertex  $v\in R$ with $\phi(v)=i$.
Second, let $X=\{x_1,x_2,x_{3}\}$ be a set of new vertices disjoint from $V(G)$.
Now, let $Y=Y(G,R,\phi)$ be the graph with vertex set $ (V( G) - R) \cup X$,
such that $Y[V(G)-R]=G-R$ and
$N(x_i) = R'_i \cup (X-x_i)$ for $i=1,2,3$. 
\end{defn}

\begin{claim} \label{coloring Y}
Suppose $R \subset V(G)$,  and $\phi$ is a $3$-coloring  of $G[R]$.
Then $\chi(Y(G, R, \phi)) \geq 4$.
\end{claim}
{\bf Proof.}
Let $G' = Y(G, R, \phi)$.
Suppose $G'$ has a $3$-coloring $\phi':V(G') \rightarrow C=\{1,2,3\}$.
By construction of $G'$, the colors of all $x_i$ in $\phi'$ are distinct.
So we may assume that $\phi'(x_i) = i$ for $1 \leq i \leq 3$.
By construction of $G'$, for all vertices $u \in R'_i$, $\phi'(u) \neq i$.
Therefore $\phi|_R \cup \phi'|_{V(G)-R}$ is a proper coloring of $G$,  a contradiction.
\qed

\begin{claim} \label{very small potential}
There is no $R \subsetneq V(G)$ with $|R| \geq 2$ and $\rho_{G}(R) \leq 5$.
\end{claim}
{\bf Proof.} Let $2\leq |R|<|V(G)|$ and
 $\rho(R)=m=\min\{\rho(W)\,:\, W\subsetneq V(G),\; |W|\geq 2\}$.
Suppose $m\leq 5$. Then $|R| \geq 4$. Since $G$ is $4$-critical,
 $G[R]$ has a proper coloring $\phi:R \rightarrow C = \{1, 2,3\}$.  
Let $G'=Y(G,R,\phi)$. By Claim~\ref{coloring Y},  $G'$ is not $3$-colorable.
Then it contains a $4$-critical subgraph $G''$. Let $W=V(G'')$. Since $|R| \geq 4>|X|$,
 $|V(G'')|<|V(G)|$.
So, by the minimality
of $G$, $\rho_{G''}(W) = \rho_{G'}(W)\leq 2$.
Since $G$ is $4$-critical by itself, $W\cap X\neq \emptyset$.
Since every non-empty subset of $X$ has potential at least $5$,
  $\quad\rho_{G}((W - X) \cup R) \leq \rho_{G'}(W)-5+m\leq m-3.$
Since $(W-X)\cup R\supset R$, $|(W-X)|\cup R|\geq 2$. Since $\rho_{G}((W - X) \cup R)<\rho_{G}( R)$,
by the choice of $R$, $(W-X)\cup R=V(G)$. But then $\rho_{G}(V(G)) \leq m-3\leq 2$,
 a contradiction. \qed

\begin{claim} \label{small potential}
If $R \subsetneq V(G)$, $|R|\geq 2$ and $\rho(R) \leq 6$, then $R$ is a $K_{3}$.
\end{claim}
{\bf Proof.}
Let $R$ have the smallest $\rho(R)$ among $R \subsetneq V(G)$, $|R|\geq 2$.
Suppose $m=\rho(R) \leq 6$ and $G[R]\neq K_{3}$. Then $|R|\geq 4$.
By Claim~\ref{very small potential}, $m=6$.

Let $R_* =\{u_1,\ldots,u_s\}$ be the set of vertices in $R$ that have neighbors outside of  $R$.
Because 
$G$ is $2$-connected, $s \geq 2$. Let  $H=G[R]+u_1u_2$.
Since $R\neq V(G)$, $|V(H)|<|V(G)|$.
By the minimality of $\rho(R)$, for every $U\subseteq R$ with $|U|\geq 2$,
$\rho_{H}(U)\geq \rho_{G}(U)-3\geq \rho_{G}(R)-3 \geq 3.$
Thus $P(H)\geq 3$, and by (\ref{j26}), $H$ has
 a proper $3$-coloring $\phi$ with colors in $C= \{1, 2,3\}$.
Let $G'=Y(G,R,\phi)$.
Since $|R|\geq 4$, $|V(G')|<|V(G)|$. 
By Claim \ref{coloring Y}, $G'$ is not  $3$-colorable. Thus
$G'$ contains  a $4$-critical subgraph $G''$. Let $W=V(G'')$. By the minimality
of $|V(G)|$, $\rho_{G''}(W) = \rho_{G'}(W)\leq 2$.
Since $G$ is $4$-critical by itself, $W\cap X\neq \emptyset$.
By Fact~\ref{list of potentials},
 if $|W \cap X| \geq 2$ then  $\rho_{G}((W - X) \cup R) \leq \rho_{G'}(W)-6+6\leq 2$,  a contradiction again.
So, we may assume that $X \cap W = \{x_1\}$.
Then
\begin{equation}\label{e2 - k}
\rho_{G}((W - \{x_1\}) \cup R)  \leq  (\rho_{G'}(W)-5) + \rho_{G}(R)\leq    \rho_{G}(R) - 3.
\end{equation}
By the minimality of $\rho_{G}(R)$, $(W - \{x_1\})\cup  R = V(G)$.
This implies that $W = V(G') - X+x_1$.

Let $R_1 = \{u \in R_*: \phi(u) = \phi(x_1)\}$.
If $|R_1|=1$, then
$\rho_{G}(W-x_1 \cup R_1)= \rho_{H}(W) \leq 2,$
a contradiction.
Thus, $|R_1|\geq 2$.
Since $R_1$ is an independent set in $H$ 
and $u_1u_2\in E(H)$,
 we may assume that $u_2\notin R_1$.
Then an  edge $u_2z$ connecting $u_2$ with $V(G)-R$
was not accounted in (\ref{e2 - k}).
So, in this case  instead of (\ref{e2 - k}), we have
$$
	\rho_{G}((W - \{x_1\}) \cup R)  \leq  \rho_{G'}(W)-5-3 + \rho_{G}(R)
	  \leq \rho_{G}(R)-6\leq 0.\qed$$

\begin{claim}\label{cliques}
$G$ does not contain $K_4-e$.
\end{claim}
{\bf Proof.} If $G[R]=K_4-e$, then $\rho_G(R)=5(4)-3(5)=5$, a contradiction to Claim~\ref{small potential}.\qed

\begin{claim} \label{clusters in cliques}
Each triangle in $G$ contains at most one vertex of degree $3$.
\end{claim}
{\bf Proof.}
By contradiction, assume that $G[\{x_1,x_2,x_3\}]=K_3$ and
$d(x_1) = d(x_2) = 3$. Let
 $N(x_1) =X-x_1+a$ and  $N(x_2) =X-x_2+b$.
By Claim~\ref{cliques},   $a \neq b$.
 Define $G' = G - \{x_1,x_2\} + ab$.
Because $\rho_{G}(W) \geq 6$ for all $W \subseteq G-\{x_1,x_2\}$ with $|W| \geq 2$, and adding an edge decreases the potential of a set by $3$,
$P(G') \geq \min\{(5,6-3\}=3.$
So, by (\ref{j26}), $G'$
 has a proper $3$-coloring $\phi'$ with $\phi'(a) \neq \phi'(b)$.
This easily extends to a proper $3$-coloring of $V(G)$.
\qed

\begin{claim} \label{adjacent k-1}
Let $xy \in E(G)$ and $d(x)=d(y)=3$.
Then both, $x$ and $y$ are in 
 triangles.
\end{claim}
{\bf Proof.}
Assume that $x$ is not in a $K_3$.
Suppose $N(x)=\{y,u,v\}$. Then $uv\notin E(G)$.
Let $G'$ be obtained from $G-y-x$ by gluing $u$ and $v$ into a new vertex $u*v$.
Since $|V(G')|<|V(G)|$, $G'$ is smaller than $G$.
If  $G'$ has a $3$-coloring $\phi':V(G') \rightarrow C = \{1, 2, 3\}$, then
we extend it to a proper $3$-coloring $\phi$ of $G$ as follows:
define $\phi|_{V(G)-x-y-u-v} = \phi'|_{V(G')-u*v}$, then
let $\phi(u)=\phi(v)=\phi'(u*v)$,
  choose $\phi(y) \in C - (\phi'(N(y) - x))$,
and $\phi(x) \in C-\{\phi(y), \phi(u)\} $.

So, $\chi(G')\geq 4$ and $G'$ contains a $4$-critical subgraph $G''$. Let $W=V(G'')$.
Since $G''$ is smaller than $G$, $\rho_{ G''}(W) = \rho_{ G'}(W) \leq 2$. Since $G''$ is not a subgraph of $G$,
 $u*v\in W$. Let $W'=W-u*v+u+v+x$.
Then
$\quad\rho_{G}(W') \leq 2 +5(2) -3(2) = 6,\quad$
since $G[W']$ has two extra vertices and at least two extra edges in comparison
with $G''$.
This contradicts Claim \ref{small potential} because $y \notin W'$ and so $W'\neq V(G)$.
\qed

By Claims \ref{cliques} and \ref{adjacent k-1}, we have
\begin{equation}\label{41}
\mbox{Each vertex with degree $3$ has at most $1$ neighbor with degree $3$.}
\end{equation}



We will now use discharging to show that $|E(G)| \geq \frac{5}3 |V(G)|$, which will finish the proof 
 of Theorem \ref{4-critical}.
Each vertex begins with charge equal to its degree.
If $d(v) \geq 4$, then $v$ gives charge $\frac16$ to each neighbor with degree $3$.
Note that $v$ will be left with charge at least $\frac56 d(v) \geq \frac{10}{3}$.
By 
(\ref{41}), each vertex of degree $3$ will end with charge at least $3 + \frac26 = \frac{10}3$.
\qed

\section{Proof of Theorem~\ref{gr}}
 Let $G$ be a plane graph with fewest elements (vertices and edges) for which
the theorem does not hold. Then $G$ is $4$-critical and in particular $2$-connected.
Suppose $G$ has $n$ vertices, $e$ edges and $f$ faces.

CASE 1: $G$ has no $4$-faces. Then $5f\leq 2e$ and so $f\leq 2e/5$. By
this and Euler's Formula $n-e+f=2$, we have $n-3e/5\geq 2$, i.e., $e\leq \frac{5n-10}{3}$,
a contradiction to Theorem~\ref{k-critical}.

CASE 2: $G$ has a $4$-face $(x,y,z,u)$. Since $G$ has no triangles, $xz,yu\notin E(G)$.
 If the graph $G_{xz}$ obtained from $G$ by gluing $x$ with $z$
 has no triangles, then by the minimality of $G$, it is $3$-colorable, and
so $G$ also is $3$-colorable. Thus $G$ has an $x,z$-path $(x,v,w,z)$ of length $3$. Since $G$ itself
has no triangles, $\{y,u\}\cap \{v,w\}=\emptyset$ and there are no edges between $\{y,u\}$ and $\{v,w\}$.
But then $G$ has no $y,u$-path  of length $3$, since such a path must cross the path $(x,v,w,z)$.
Thus the graph $G_{yu}$ obtained from $G$ by gluing $y$ with $u$
 has no triangles, and so, by the minimality of $G$, is $3$-colorable.
Then  $G$ also is $3$-colorable, a contradiction.\qed

\medskip
{\bf Acknowledgment.} We thank Michael Stiebitz for helpful comments.

\end{document}